\documentclass[11pt,leqno]{article}
 \usepackage{amssymb,amsfonts}
 \usepackage{amsmath,amsthm,amsxtra}
 \usepackage[all]{xy}

 \newcommand{\QQ}{\mathbb{Q}}
 \newcommand{\ZZ}{\mathbb{Z}}

 \newcommand{\F}{\mathcal{F}}
 \newcommand{\J}{\mathcal{J}}

 \newcommand{\Fract}{\mathop{\rm Fract \, }}
 \newcommand{\ord}{\mathop{\rm ord \, }}

 \makeatletter
 \renewcommand\section{\@startsection {section}{1}{\z@}%
                                    {-3.5ex \@plus -1ex \@minus -.2ex}%
                                    {2.3ex \@plus.2ex}%
                                    {\normalfont\large\bfseries}}
 \renewcommand\subsection{\@startsection{subsection}{2}{\z@}%
                                      {-3.25ex\@plus -1ex \@minus -.2ex}%
                                      {0ex \@plus .0ex}%
                                      {\normalfont\normalsize\bfseries}}

 \newtheorem{theorem}{Theorem}
 \newtheorem{lemma}{Lemma}

 \theoremstyle{remark}

 \newtheorem{corollary}{Corollary}

 \newcommand{\alphaparenlist}{
   \renewcommand{\theenumi}{\alph{enumi}}%
   \renewcommand{\labelenumi}{(\theenumi)}%
 }

 \makeatletter
 \def\@maketitle{\newpage
  \null
  \vskip 2em
  \begin{center}%
   {\@date}%
   \vskip 3em
   {\Large\bf \@title \par}%
   \vskip 1.5em
   {\normalsize
    \lineskip .5em
    \begin{tabular}[t]{c}\@author
    \end{tabular}\par}%
   \vskip 2em

  \end{center}%
  \par
  \vskip 2.5em}
 \makeatother

\begin{document}

 \title{A differential analog of a theorem of Chevalley}
 \author{Victor G. Kac\thanks{Department of Mathematics, M.I.T.,
     Cambridge, MA 02139, USA.~~kac@math.mit.edu}~~\thanks{Supported in part by NSF grant
     DMS-9970007.}}

 \maketitle

\begin{abstract}
  In this note a proof of a differential analog
 of Chevalley's theorem \cite{C} on homomorphism extensions is given.
 An immediate corollary is a condition of finitenes of extensions
 of differential algebras and several equivalent definitions of
 a differentially closed field, including Kolchin's Nullstellensatz.

\end{abstract}

 In this note I give a proof of the following differential analog
 of Chevalley's theorem \cite{C} on homomorphism extensions.

 \begin{theorem}
   \label{th:1}
   Let $S$ be a differential algebra over $\QQ$ with no zero
   divisors and let $b$ be a non-zero element of $S$.  Let $R$ be
   a differential subalgebra of $S$ over which $S$ is
   differentially finitely generated.  Let $\F$ be a
   differentially closed field of characteristic $0$.  Then there
   exists a non-zero element $a$ of $R$ such that any homomorphism
   $\varphi : R \to \F$ which does not annihilate $a$ extends to a
   homomorphism $\psi : S \to \F$ which does not annihilate $b$.
 \end{theorem}

 An almost immediate consequence of the proof of
 Theorem~\ref{th:1} is

 \begin{theorem}
   \label{th:2}
   Let $\F$ be a differentially closed field of characteristic $0$
   and let $S \supset R$ be differentially finitely generated
   differential algebra and subalgebra over $\F$.  Suppose that
   there exists a non-zero element $b$ of $S$ such that any
   homomorphism $\varphi :R \to \F$ has only finitely many
   extensions $\psi : S \to \F$
   satisfying $\psi (b) \neq 0$.  Then the field extension $\Fract
   S \supset \Fract R$ is finite.
   In particular, if any homomorphism $\varphi : R \to \F$ has at
   most $d$ extensions $\psi : S \to \F$ with $\psi (b) \neq 0$,
   then the degree of $\Fract S $ over $\Fract R$ is at most $d$.
 \end{theorem}

 An immediate corollary of Theorem~\ref{th:1} is Kolchin's
 Nullstellensatz \cite{K} and its earlier weaker versions by Ritt
 \cite{R}, Cohn \cite {Cohn} and Seidenberg \cite{S}.

 As far as I can understand it, Theorem~\ref{th:1} is closely
 related to Blum's elimination of quantifiers theorem \cite{Blum}, \cite{M}
 in the model theory of differentially closed fields.

 In Section~1 I explain the necessary background on
 Differential Algebra, in Sections~2 and~3
 give proofs of Theorems~\ref{th:1} and~\ref{th:2} and in
 Section~4 give about a dozen of equivalent definitions
 of a differentially closed field of characteristic $0$
 (Theorem~\ref{th:1} and Kolchin's Nullstellensatz being among
 them).

 This note is an offshoot of a course in Differential Algebra that
 I gave at M.I.T. in the fall of 2000.  The motivation for teaching
 this course came from close connections of the theory of
 conformal algebras to differential algebras \cite{BDK}.  Namely,
 each Lie conformal algebra defines a functor from the category of
 differential algebras to the category of Lie algebras
 (very much like a group scheme defines a
 functor from the category of algebras to the category of
 groups).  I am grateful to my students for their enthusiasm,
 especially to A.~De~Sole for several corrections and improvements.
 I would like to thank A.~Buium for very useful and enlightening
 correspondence and J.~Young for bringing Marker's paper \cite{M}
 to my attention and explaining some parts of it.

 \vspace{3ex}

 \textbf{1.~~}
 Here I recall some terminology and facts from Differential
 Algebra.  All of this can be easily found in two excellent books
 \cite{K} and \cite{B}, the primary source of both being Ritt's
 foundational book \cite{R}.

 By an algebra we always mean a commutative associative unital
 algebra.  A \emph{differential algebra} $R$ is an algebra over a
 field with a fixed derivation $\delta$.  By a homomorphism of
 differential algebras we always mean a homomorphism commuting
 with derivations.  Likewise subalgebras and ideals are assumed
 to be $\delta$-invariant, though to emphasize this we shall often
 call them differential subalgebras and ideals.

 The first important observation of Ritt is that in a differential
 algebra over $\QQ$ the radical of a differential ideal is a
 differential ideal (see \cite{R}, \cite{Ra}, \cite{K}
 Lemma~1.8).  (A counterexample in characteristic $p$ is the zero
 ideal in $R=\F[x]/(x^p)$, $\delta = d/dx$.)  Another
 important fact is the differential Krull theorem:  any
 differential radical ideal is an intersection of differential
 prime ideals (see \cite{Ra}, \cite{K} Theorem~2.1).  These facts
 imply, in particular, that maximal (among) differential ideals are
 prime.

 If a differential algebra $R$ has no zero divisors, we can form
 the field of fractions $\Fract R$ and extend the derivation
 $\delta$ to this field.

 One says that a differential algebra $S$ is differentially
 \emph{generated} by elements $x_1 , \ldots , x_n$ over a
 differential subalgebra $R$, and writes $S=R \{ x_1 , \ldots ,
 x_n \}$, if the algebra $S$ is generated by all elements from $R$
 and all derivatives $x^{(k)}_i$, $i=1,\ldots ,n$, $k \in \ZZ_+
 = \{ 0,1,2,\ldots \}$, of the elements~$x_i$.

 Let $S = R \{ x \}$ be a differential algebra with no zero
 divisors.  One says that $x$ is \emph{differentially
   transcendental} over $R$ if all elements $x^{(k)}$, $k \in
 \ZZ_+$, are algebraically independent over $\Fract R$; otherwise
 $x$ is called \emph{differentially algebraic} over $R$.

 The algebra of \emph{differential polynomials} over a
 differential algebra $R$ in the differential indeterminates
 $y_1,\ldots ,y_n$ is the algebra of polynomials $R [y^{(k)}_i,
 i=1,\ldots ,n, \,\, k \in \ZZ_+]$ with the derivation $\delta$
 extended from $R$ by the rule $\delta (y^{(k)}_i) = y^{(k+1)}$.
 This algebra is denoted by $R \{ y_1 ,\ldots ,y_n \}$.

 Consider the differential algebra $R \{ y \}$ of differential
 polynomials over $R$ in one differential indeterminate $y$.
 Let $A(y) \in R \{ y \} \backslash R$, be a ``non-constant''
 differential polynomial.  The largest $r$ for which $y^{(r)}$ is
 present in $A(y)$ is called the \emph{order} of $A(y)$ and is
 denoted by $\ord A$.  One can write in a unique way:
 \begin{displaymath}
   A(y) = I_A (y) y^{(r)d} + I_1 (y) y^{(r)d-1}
   + \cdots + I_d (y) \, ,
 \end{displaymath}
 where $\ord I_A (y) < r$, $\ord I_j(y)<r$ and $I_A (y) \neq 0$.
 Then $d$ is called the \emph{degree} of $A(y)$ and is denoted by
 $\deg A(y)$, and $I_A(y)$ is called the \emph{initial} of
 $A(y)$.  The differential polynomial $S_A (y) = \partial A(y)
 /\partial y^{(r)}$ is called the \emph{separant} of $A(y)$.
 The important property of the characteristic 0 case is that
 $S_A(y) \neq 0$ if $A \notin R$.

 For $A,B \in R \{ y \} \backslash R$ we write $A<B$ if either
 $\ord A < \ord B$, or $\ord A =\ord B$ and $\deg A < \deg B$.  We
 also write $A<B$ if $A \in R$, $B \notin R$.  Note that $I_A<A$
 and $S_A<A$.

 The basic result of Ritt's theory \cite{R} (see \cite{K}
 Lemma~7.3 or \cite{B} (2.3)) is the following division
 algorithm:

 Given a non-constant differential polynomial $A (y)$, for any
 $F(y) \in R \{ y \}$ there exists\break $G(y) \in R \{ y \}$ with $G<A$
 and $m,n \in \ZZ_+$ such that:
 \begin{displaymath}
   I_A (y)^m S_A (y)^n F(y) \equiv
    G(y) \mod [A(y)] \, ,
 \end{displaymath}
 where $[A(y)]$ is the ideal of $R \{ y \}$ generated by
 $A(y)^{(k)}$, $k \in \ZZ_+$.  Furthermore, one can take $m=0$ if
 only $\ord G \leq \ord A$ is required, which is called the weak
 division algorithm.

 The important notion of a differentially closed field was
 introduced by Robinson \cite{Rob}.  His axioms have been
 considerably simplified by Blum \cite{Blum}, and it is her
 definition, given below, that is commonly used.  A differential
 field $\F$ is called \emph{differentially closed} if for any two
 differential polynomials $A(y), B(y) \in \F \{ y \}$ with $B \neq
 0$ and $\ord B <\ord A$ there exists $\alpha \in \F$ such that $A
 (\alpha)=0$ and $B (\alpha) \neq 0$.

 Any differentially closed field is algebraically closed
 (i.e.,~has no non-trivial finite extensions), but it always has
 differentially algebraic extensions.  Nevertheless, it turned out
 to be the right substitute for Differential Algebra of the notion
 of an algebraically closed field.

 The existence of a differentially closed field containing a given
 differential field of characteristic 0 is easy to establish in the framework
 of Model Theory (see e.g.~\cite{M}).  An elementary proof (i.e.,~without
 reference to model theory) may be found in \cite{B} (5.2).

 \vspace{3ex}

 \textbf{2.~~}\emph{Proof of Theorem~1.~~}
 The general scheme of the proof is the same as Chevalley's
 \cite{C}.  We have:
 \begin{equation}
   \label{eq:1}
   S=R \{ x_1 ,\ldots ,x_{n-1} \} \{ x_n \} \supset
   R \{ x_1 ,\ldots , x_{n-1} \} \supset R \, .
 \end{equation}
 By induction on $n$ we reduce the proof to the case $n=1$.
 Indeed, from being true for $n=1$, we conclude that there exists
 a non-zero $b_1 \in R \{ x_1 , \ldots , x_{n-1} \}$ such that any
 homomorphism $\psi_1 : R \{ x_1 , \ldots , x_{n-1}\} \to \F$ with
 $\psi_1 (b_1) \neq 0$ extends to $\psi : R \{ x_1 , \ldots ,
 x_n \} \to \F$ with $\psi (b) \neq 0$, and by the inductive
 assumption, there exists a non-zero $a \in R$ such that any
 homomorphism $\varphi : R \to \F$ with $\varphi (a) \neq 0$
 extends to $\psi_1 :R \{ x_1 ,\ldots ,x_{n-1} \} \to \F$ with
 $\psi_1 (b_1) \neq 0$.

 Thus, we may assume that $S=R \{ x\}$.  Given a homomorphism
   $\varphi :R \to \F$, we may extend it to the algebras of
   differential polynomials $A(y) \mapsto A^{\varphi}(y)$ by
   applying $\varphi$ to coefficients.  Let $B(y) \in R \{ y \}$
     be such that $B(x) =b$.  We consider separately two cases.

 \textbf{Case 1:}~~$x$ is differentially transcendental over $R$.  Let $a \in
       R$ be any non-zero coefficient of $B(y)$.  If $\varphi (a)
       \neq 0$, then $B^{\varphi}(y) \in \F \{ y \}$ is a non-zero
       differential polynomial, hence there exists $\alpha \in \F$
       which is not a root of $B^{\varphi}(y)$ (we take
       $B=B^{\varphi}(y)$ and $A$ with $\ord A > \ord B$ in the
       definition of a differentially closed field).  Then $\psi
       (Q(x)) = Q^{\varphi} (\alpha)$ is a well defined
       homomorphism $S \to \F$ with $\psi (b) \neq 0$.

 \textbf{Case 2:}~~$x$ is differentially algebraic over $R$.
 Let $A (y) \in R \{ y
 \}$ be a minimal in the partial ordering $<$ irreducible (in the
 usual sense) over $\Fract R$ differential polynomial such that
 $A(x) =0$.  If $F (y) \in R \{ y \}$ is such that $F(x)=0$, apply
 the division algorithm:
 \begin{displaymath}
   S_A (y)^m I_A (y)^n F(y) \equiv
   G(y) \mod [A(y)] \, ,
 \end{displaymath}
 where $G<A$.  Since $F(x) =0$ and $A(x) = A'(x) = \cdots =0$, we
 see that $G (x) =0$, hence, due to minimality of $A$, $G(y)=0$,
 and we have
 \begin{equation}
   \label{eq:2}
   S_A (y)^m I_A (y)^n F(y) \in [A(y)] \, .
 \end{equation}

 Let $a_1 \in R$ be a non-zero coefficient of $I_A(y)$.  Let $D
 (y) \in R \{ y \}$ be the discriminant of $A(y)$ viewed as a
 polynomial in $y^{(\ord A)}$.  Note that $D (y) \neq 0$ since
 $A(y)$ is an irreducible polynomial, and that $\ord D(y) <\ord
 A(y)$.  Let $a_2 \in R$ be a non-zero coefficient of $D(y)$.

 Suppose that $\varphi (a_1a_2) \neq 0$.  Then $\ord A^{\varphi}
 (y) = \ord A (y) > \ord I^{\varphi}_A (y) D^{\varphi}(y)$ and\break
 $I^{\varphi}_A (y) D^{\varphi}(y) \neq 0$.  Since $\F$ is
 differentially closed, there exists $\alpha \in \F$ which is a
 root of $A^{\varphi}(y)$ but not a root of $I^{\varphi}_A (y)
 D^{\varphi}(y)$.  Since $D^{\varphi} (\alpha) \neq 0$,
 $A^{\varphi} (y)$ and $S^{\varphi}(y)$ have no common roots,
 hence $S^{\varphi}_A (\alpha) \neq 0$.  Since also $I^{\varphi}_A
 (\alpha) \neq 0$, but $A^{\varphi}_A (\alpha)=0$, we conclude
 from (\ref{eq:2}) that $F^{\varphi} (\alpha)=0$.  Therefore $\psi
 (Q(x))=Q^{\varphi}(\alpha)$ is a well defined homomorphism $S \to
 \F$ which extends $\varphi :R \to \F$.

 It remains to take care of the condition $\psi (b) \neq 0$ by an
 appropriate choice of $\alpha$ (satisfying the above conditions
 as well).  By the weak division algorithm we have:
 \begin{equation}
   \label{eq:3}
   S^n_A (y) B(y) = B_1(y) \mod [A (y)] \, ,
 \end{equation}
 where $n \in \ZZ_+$ and $\ord B_1(y) \leq \ord A(y)$.  Letting
 $y=x$ in (\ref{eq:3}) we get
 \begin{equation}
   \label{eq:4}
   B_1(x) = S^n_A (x) b \, .
 \end{equation}
 Since $S^{\varphi}_A (\alpha) \neq 0$, we conclude that $S_A (x)
 \neq 0$, hence $B_1 (x) \neq 0$, due to (\ref{eq:3}).  Let $r(y)
 \in R \{ y \}$ be the resultant of the polynomials $B_1(y)$ and $
 A(y)$ viewed as polynomials in $y^{(\ord A)}$.  Since $A(x) =0$
 and $B_1 (x) \neq 0$, we conclude that $r(y) \neq 0$ (otherwise $
 A(y)$ and $B_1(y)$ would have a common root in $\Fract R$ and
 therefore $B_1(y)$ would be divisible by $A(y)$ due to its
 irreducibility).  Note that $\ord r(y) < \ord A(y)$.  Let $a_3
 \in R$ be a non-zero coefficient of $r(y)$.

 We let $a=a_1a_2a_3$ and suppose that $\varphi (a) \neq 0$.
 Choose $\alpha \in \F$ such that, as before, $A^{\varphi}(\alpha)
 \neq 0$, $I^{\varphi}_A (\alpha) D^{\varphi} (\alpha) \neq 0$,
 and, in addition, $r^{\varphi} (\alpha) \neq 0$.  As before,
 define $\psi (Q(x))=Q^{\varphi}(\alpha)$.  As before, this is a
 well defined homomorphism $S \to \F$, and, by (\ref{eq:4}):
 $S^{\varphi}_A (\alpha) \psi (b) = B^{\varphi}_1 (\alpha)$.  But
 $B^{\varphi}_1 (\alpha) \neq 0$, since $r^{\varphi}(\alpha)\neq
 0$ and therefore $A^{\varphi}(y)$ and $B^{\varphi}_1 (y)$ have no
 common roots.  Hence $\psi (b) \neq 0$.

     \begin{corollary}
       \label{cor:1}
 If $S$ is a differentially finitely generated differential
 algebra over $\F$, then for any non-zero $b \in S$ there exists
 an $\F$-algebra homomorphism $\psi : S \to \F$ such that $\psi
 (b) \neq 0$.
     \end{corollary}

 \vspace{3ex}
 \textbf{3.~~}
 In order to prove Theorem~\ref{th:2}, we need the following
 lemma.

 \begin{lemma}
   \label{lem:1}
 Let $\F$ be a differentially closed filed.  Then

 \alphaparenlist
 \begin{enumerate}
 \item 
 For any non-zero $A(y) \in \F \{ y \}$ there exists infinitely
 many $\alpha \in \F$ such that $A(\alpha) \neq 0$.

 \item  
   If $A(y) \in \F \{ y \} \backslash \F$ and there exists $B(y)
   \in \F \{ y \}$ such that $B(y) \neq 0$, and $B(y) < \ord A
   (y)$ and only finitely many roots of $A(y)$ are not roots of
   $B(y)$, then $\ord A(y) =0$.
 \end{enumerate}
 \end{lemma}

 \emph{Proof.~~}
 (a) Take a sequence $A_j(y) \in \F \{ y \}$, $j \in \ZZ_+$, of
 increasing order, such that $A_0 (y) = A(y)$.  For each
 $j=1,2,\ldots$ there exists $\alpha_j \in \F$ which is a root
 of $A_j$, but not of $A_0 \ldots A_{j-1}$.  Hence all
 $\alpha_1,\alpha_2,\ldots$ are not roots of $A(y)$.

 (b) Let $\alpha_1,\ldots , \alpha_n \in \F$ be all roots of $A(y)$
 which are not roots of $B(y)$.  Note that $n \geq 1$ since $\F$ is differentially closed.  Let
 $\displaystyle{B_1 (y) = \prod^n_{i=1} (y-\alpha_i) B(y)}$ and suppose that
 $\ord A(y) >0$.  Then $\ord B_1 (y) < \ord A(y)$, hence there
 exists $\beta \in \F$ which is a root of $A(y)$, but not of $B_1
 (y)$.  But then $\beta \neq \alpha_i$ for all $i$, a
 contradiction.
 \hfill $\square$
 \vspace{2ex}

 \emph{Proof of Theorem~2.~~}
 Again, using (\ref{eq:1}), we reduce the proof to the case $S=R
 \{ x \}$.  Certainly any homomorphism $R \{ x_1 , \ldots ,
   x_{n-1} \} \to \F$ extends in only finitely many ways to $S \to
   \F$ that does not annihilate $b$.  Hence from $n=1$ case we
   conclude that $\Fract S$ is finite over $\Fract R \{ x_1 ,
   \ldots , x_{n-1} \}$.  By Theorem~\ref{th:1}, there exists $b'
   \neq 0$ in $R \{ x_1 , \ldots , x_{n-1} \}$ such that any
   homomorphism $\varphi$ to $\F$ with $\varphi (b') \neq 0$
   extends to $\psi : S \to \F$ with $\psi (b) \neq 0$.  Hence
   there exists only finitely many homomorphisms $\psi : R \{ x_1
   , \ldots , x_{n-1} \} \to \F$ which extend $\varphi :R \to \F$
   such that $\psi (b')\neq 0$.  Hence, by the inductive
   assumption, $\Fract R \{ x_1 , \ldots , x_{n-1} \}$ is finite
   over $\Fract R$ and $\Fract S$ is finite over $\Fract R$.

 Since $R$ is differentially finitely generated over $\F$, by
 Corollary~\ref{cor:1}, for any non-zero $a \in R$ we can find a
 homomorphism $\varphi : R \to \F$ with $\varphi (a) \neq
 0$. Again we consider separately two cases of $S = R \{ x \}$,
 keeping notations of the proof of Theorem~\ref{th:1}.

 \textbf{Case 1:}  $x$ is differentially transcendental over $R$.
 Let $a \in R$ be a non-zero coefficient of $B(y)$, so that
 $B^{\varphi} (y) \neq 0$.  By Lemma~\ref{lem:1}a,
 $B^{\varphi}(y)$ has infinitely many non-roots $\alpha$, and for
 each of them the homomorphism $\psi (B(x)) =B^{\varphi}(\alpha)$
 extends $\varphi$ such that $\psi (b) \neq 0$.  Thus, this case
 is impossible.

 \textbf{Case 2:}  $x$ is differentially algebraic over $R$.  Take
 $a=a_1a_2a_3$ from the proof of Theorem~\ref{th:1} and $\varphi
 :R \to \F$ with $\varphi (a) \neq 0$.  From the proof of
 Theorem~\ref{th:1} we know that $\varphi$ extends to $\psi : S
 \to \F$ with $\psi (b) \neq 0$ by letting $\psi (Q(x))=Q^{\varphi}
   (\alpha)$, where $\alpha$ is a root of $A^{\varphi} (y)$ and
 is not a root of $I^{\varphi}_A (y) D^{\varphi}(y)
 r^{\varphi}(y)$.  Since, by conditions of the theorem, there
 exists only finitely many such $\alpha$, we conclude by
 Lemma~\ref{lem:1}, that $\ord A^{\varphi}(y)=0= \ord A(y)$.  In
 other words, $A(y)$ is an ordinary (non-zero) polynomial over $R$
 such that $A(x)=0$.  Hence $\Fract S$ is finite over $\Fract R$.

 \hfill $\square$

 \vspace{3ex}

 \textbf{4.~~Theorem~3.~~}
 The following properties of a differential filed $\F$ of
 characteristic $0$ are equivalent:

 \alphaparenlist

 \begin{enumerate}
 \item 
   $F $ is differentially closed.

 \item 
   If $A,B \in \F \{ y \}$ are such that $A$ is irreducible, $B$
   is not divisible by $A$ and $\ord A \geq \ord B$, then there
   exists $\alpha \in \F$ such that $A (\alpha) =0$, $B
   (\alpha)\neq 0$ and $S_A (\alpha) \neq 0$.

 \item 
   If $\J$ is a prime differential ideal of $\F \{ y \}$ and $B \in
   \F \{ y \} \backslash \J$, then there exists $\alpha \in \F$
   such that $f (\alpha)=0$ for all $f \in \J$ and $B (\alpha)
   \neq 0$.

 \item 
 The same as (c), where $\J$ is a radical differential ideal.

 \item 
   The same as (c) (resp. d) where $\F \{ y \}$ is replaced by $\F
   \{ y_1 , \ldots , y_n \}$, $n \geq 1$.

 \item 
   If $\J$ is a proper prime differential ideal of $\F \{
   y_1,\ldots ,y_n \}$, $n \geq 1$ arbitrary, then there exists
   $\alpha \in \F$ such that $f (\alpha)=0$ for all $f \in \J$.

 \item 
   The same as (f), where $\J$ is a proper differential ideal
   (resp. radical differential ideal).

 \item 
   The same as (f), where $\J$ is a proper maximal differential
   ideal.

 \item 
   Given any $S$, $R$ and $b$ as in Theorem~\ref{th:1}, the
   conclusion of Theorem~\ref{th:1} holds for $\F$.
 \end{enumerate}

 \emph{Proof.~~}
   (a) implies (i) by Theorem~\ref{th:1}.  (i) implies (c) by
   taking $S=\F \{ y \} /\J$, $b=$ image of $B$ in $S$ and
   applying Corollary~\ref{cor:1} (and the differential Krull
   theorem). (e)~trivially implies~(c) and~(d).  A standard
   argument of adding extra indeterminates shows that (f) is
   equivalent to~(e).  By the differential Krull theorem, (g)~is
   equivalent to~(f).  Since every prime differential ideal can be
   included in a maximal one, which is prime too, (h)~is
   equivalent to~(f).  (b)~trivially implies~(a) (by taking in (a)
   an irreducible factor of $A$ of the same order).  Finally, the
   implication (c) $\Rightarrow$ (b) is proved in the same way as
   3)~$\Rightarrow$~1) in Theorem~5.1 of \cite{B}.  Indeed, $\{ A
   \} : S_A$ is a prime ideal by Ritt's structure theorem \cite{B}
   Theorem~2.5 or \cite{K} Lemma~7.9, and this ideal does not
   contain $S_A B$ by Ritt's divisibility theorem \cite{B}
   Theorem~2.4 or \cite{K} Lemma~7.8.
 \hfill $\square$
 \vspace{2ex}

 \emph{Remark.~~}All proofs can be extended without difficulty to
 the case of several commuting derivations $\delta_1,\ldots ,\delta_m$.

 \end{document}